\documentclass[11pt]{article}
\usepackage{amsmath}
\usepackage{amsfonts}
\usepackage{amssymb}
\usepackage{hyperref}
\usepackage{indentfirst}
\usepackage{mathrsfs}

\usepackage[top=1in,bottom=1in,left=1.25in,right=1.25in]{geometry}
\date{}
\allowdisplaybreaks

\begin{document}

\renewcommand{\baselinestretch}{1.2}
\renewcommand{\arraystretch}{1.0}

\title{\bf  Drinfel'd twist of multiplier Hopf algebras}
\author
{
  \textbf{Tao Yang} \footnote{College of Science, Nanjing Agricultural University, Nanjing 210095, Jiangsu, CHINA.
             E-mail: tao.yang@njau.edu.cn}
%  \textbf{Xuan Zhou},
%  \textbf{Xiaoyan Zhou}
}

\maketitle

\begin{center}
\begin{minipage}{12.cm}

 {\bf Abstract } This paper generalizes the Drinfel'd twist to the multiplier Hopf algebra case.
 For a multiplier Hopf algebra $A$ with a twist $J$, we construct a new multiplier Hopf algebra $A^{J}$.
 Furthermore, if $A$ is quasitriangular, then $A^{J}$ is also.
 Finally, for a counimodular algebraic quantum group $A$, $A^{J}$ is an algebraic quantum group.
\\

\\


\end{minipage}
\end{center}
\normalsize

\section*{Introduction}
\def\theequation{\thesection.\arabic{equation}}
\setcounter{equation}{0}

 Hopf algebras in symmetric monoidal categories, arise naturally in the deformation-quantization of
 triangular solutions of the classical Yang-Baxter equations. 
 This deformation-quantization in \cite{GRZ92} was achieved by means of a Drinfel'd twist.
 
 Dinfel'd twist play a important role in Hopf algebra theory, such as \cite{GM94}. In \cite{AEG02}, 
 the authors studied the properties of Drinfel'd twist for finite-dimensional Hopf algebras.
 Given a Hopf algebra or quasitriangular Hopf algebra with a Drinfel'd twist, 
 they got another new Hopf algebra and quasitriangular Hopf algebra.
 And the authors also determine how the integral of the dual to "finite-dimensional" unimodular Hopf algebra changes under a twist.
 
 How about the infinite-dimensional case? This is a obvious question. 

 As we know, multiplier Hopf algebras can be considered as a generalization of Hopf algebras, 
 and give a nice answer to the dual of infinite-dimensional Hopf algebra. 
 Hence, in this paper, we consider the Drinfel'd twist in the multiplier Hopf algebra case and give an answer to the above question.

 The paper is organized in the following way.
 In section 1, we recall some notions which we will use in the following, such as multiplier Hopf algebras, 
 algebraic quantum groups and their duality, and pairing and actions of multiplier Hopf algebras.

 In section 2, we first give the definitions of twists of multiplier Hopf algebras, which generalizing the Drinfel'd twist of Hopf algebra.
 Let $(A, \Delta)$ be a regular multiplier Hopf algebra, then $(A^{J}, \Delta^{J})$ is also a regular multiplier Hopf algebra
 with the same algebra structure and counit as $(A, \Delta)$ and the comultiplication and antipode are given by
 $\Delta^{J}(a) = J^{-1} \Delta(a) J$, $S^{J}(a) = Q^{-1}_{J} S(a) Q_{J}$ for all $a\in A$.
 Furthermore, if $A$ is quasitriangular with quasitriangular structure $\mathcal{R}$, then $A^{J}$ is also quasitriangular and its generalized $R$-matrix given by $\mathcal{R}^{J} = J^{-1}_{21} \mathcal{R} J$.
 Finally, for a counimodular algebraic quantum group $A$, $A^{J}$ is an algebraic quantum group  
 with the elements $\varphi^{J} = u_{J} \rightharpoonup \varphi$ and $\psi^{J} = \psi \leftharpoonup u^{-1}_{J}$
 as non-zero left and right integrals on $(A^{J}, \Delta^{J})$ respectively.

\section{Preliminaries}
\def\theequation{\thesection.\arabic{equation}}
\setcounter{equation}{0}

 In this section, we will recall some of the basic notions and results in the theory of algebraic quantum groups, which will be used in this paper.

\subsection{Multiplier Hopf algebras}

 Throughout this paper, all spaces we considered are over a fixed field $K$.
 Algebras may or may not have units, but should be always non-degenerate.
 For an algebra $A$, the multiplier algebra $M(A)$ of $A$ is defined as the largest algebra with unit in which $A$ is a dense ideal
 (for more details see the appendix in \cite{V94}).

 Now, we recall the definition of a multiplier Hopf algebra (see \cite{V94} for details).
 A comultiplication on algebra $A$ is a homomorphism $\Delta: A \longrightarrow M(A \otimes A)$ such that $\Delta(a)(1 \otimes b)$ and
 $(a \otimes 1)\Delta(b)$ belong to $A\otimes A$ for all $a, b \in A$. We require $\Delta$ to be coassociative in the sense that
 \begin{eqnarray*}
 (a\otimes 1\otimes 1)(\Delta \otimes \iota)(\Delta(b)(1\otimes c))
 = (\iota \otimes \Delta)((a \otimes 1)\Delta(b))(1\otimes 1\otimes c)
 \end{eqnarray*}
 for all $a, b, c \in A$ (where $\iota$ denotes the identity map).

 A pair $(A, \Delta)$ of an algebra $A$ with non-degenerate product and a comultiplication $\Delta$ on $A$ is called
 a \emph{multiplier Hopf algebra}, if the linear map $T_{1}, T_{2}$ defined by
 \begin{eqnarray}
 T_{1}(a\otimes b)=\Delta(a)(1 \otimes b), \qquad T_{2}(a\otimes b)=(a \otimes 1)\Delta(b)
 \end{eqnarray}
 are bijective.

 The bijectivity of the above two maps is equivalent to the existence of a counit and an antipode S satisfying (and defined by)
 \begin{eqnarray}
 && (\varepsilon\otimes\iota)(\Delta(a)(1\otimes b)) = ab, \qquad  m(S\otimes\iota)(\Delta(a)(1\otimes b))=\varepsilon(a)b, \label{1.1} \\
 && (\iota\otimes\varepsilon)((a\otimes 1)\Delta(b)) = ab, \qquad  m(\iota\otimes S)((a\otimes 1)\Delta(b))=\varepsilon(b)a,\label{1.2}
 \end{eqnarray}
 where $\varepsilon:A\longrightarrow k$ is a homomorphism, $S: A\longrightarrow M(A)$ is an anti-homomorphism
 and $m$ is the multiplication map, considered as a linear map from $A\otimes A$ to $A$ and extended to
 $M(A)\otimes A$ and $A\otimes M(A)$.

 A multiplier Hopf algebra $(A, \Delta)$ is called \emph{regular} if $(A, \Delta^{cop})$ is also a multiplier Hopf algebra,
 where $\Delta^{cop}$ denotes the co-opposite comultiplication defined as $\Delta^{cop}=\tau \circ \Delta$ with $\tau$ the usual flip map
 from $A\otimes A$ to itself (and extended to $M(A\otimes A)$). In this case, $\Delta(a)(b \otimes 1), (1 \otimes a)\Delta(b) \in A \otimes A$
 for all $a, b\in A$.
 By Proposition 2.9 in \cite{V98}, multiplier Hopf algebra $(A, \Delta)$ is regular if and only if the antipode $S$ is bijective from $A$ to $A$.

 We will use the adapted Sweedler notation for multiplier Hopf algebras (see \cite{V08}), e.g., write $a_{(1)} \otimes a_{(2)}b$ for $\Delta(a)(1 \otimes b)$ and $ab_{(1)} \otimes b_{(2)}$ for $(a \otimes 1)\Delta(b)$.

\subsection{Algebraic quantum groups and their dualities}

 Assume in what follows that $(A, \Delta)$ is a regular multiplier Hopf algebra. A linear functional $\varphi$ on $A$ is called left invariant if
 $(\iota \otimes \varphi)\Delta(a) = \varphi(a)1$ in $M(A)$ for all $a \in A$. A non-zero left invariant functional $\varphi$ is called a
 \emph{left integral} on $A$. A right integral $\psi$ can be defined similarly.

 In general, left and right integrals are unique up to a scalar if they exist.
 And if a left integral $\varphi$ exists, a right integral also exists, namely $\varphi \circ S$.
 Although may be different, left and right integrals are related.
 For a left integral $\varphi$, there is a unique group-like element (modular element) $\delta \in M(A)$ such that
 $\varphi(S(a)) = \varphi(a\delta)$ for all $a \in A$.

 For an algebraic quantum group $(A, \Delta)$ with integrals, define $\widehat{A}$ as the space of linear functionals on $A$
 of the form $\varphi(\cdot a)$, where $a \in A$. Then $\widehat{A}$ can be made into a regular multiplier Hopf algebra
 with a product(resp. coproduct $\widehat{\Delta}$ on $\widehat{A}$) dual to the coproduct $\Delta$ on $A$ (resp. product of $A$).
 It is called the dual of $(A, \Delta)$.
 The various objects associated with $(\widehat{A}, \widehat{\Delta})$ are denoted as for $(A, \Delta)$ but with a hat.
 However, we use $\varepsilon$ and $S$ also for the counit and antipode on the dual.
 The dual $(\widehat{A}, \widehat{\Delta})$ also has integrals, i.e., the dual is also an algebraic quantum group.
 A right integral $\widehat{\psi}$ on $\widehat{A}$ is defined by $\widehat{\psi}(\varphi(\cdot a)) = \varepsilon(a)$ for all $a\in A$.
 Repeating the procedure, i.e., taking the dual of $(\widehat{A}, \widehat{\Delta})$,
 we can get $\widehat{\widehat{A}}\cong A$ (see Theorem (Biduality) 4.12 in \cite{V98}).

 From Definition 1.5 in \cite{D10}, an algebraic quantum group $A$ is called counimodular, if the dual multiplier Hopf algebra $\widehat{A}$ is unimodular integral, i.e., $\widehat{\delta}= 1$ in $M(A)$. For a counimodular algebraic quantum group, we have $\varphi(ab)=\varphi(bS^{2}(a))$ for
 all $a, b \in A$ (see Proposition 1.6 in \cite{D10}).

\subsection{Pairing and actions}

 By using the expression $\langle a, b \rangle$ to denote the value of a functional $b\in \widehat{A}$ in the point $a\in A$,
 we get a non-degenerate pairing of multiplier Hopf algebras in the sense of \cite{DV01}.
 This pairing gives rise to a left and a right action of $A$ on $\widehat{A}$ as follows.
 \begin{eqnarray*}
 \langle a'a, b \rangle = \langle a', a\rightharpoonup b \rangle, \qquad
 \langle aa', b \rangle = \langle a', b\leftharpoonup a \rangle
 \end{eqnarray*}
 for any $a, a' \in A$ and $b\in \widehat{A}$.

%\subsection{Group-cograded multiplier Hopf algebras}
%
% Let $(A, \Delta)$ be a multiplier Hopf algebra and $G$ be a (discrete) group with unit $e$. Assume that there is a family of (non-trivial) subalgebras $\big(A_{p}\big)_{p\in G}$ of $A$ so that
% \begin{enumerate}
% \item[(1)] $A=\bigoplus_{p\in G}A_{p}$ with $A_{p}A_{q}=0$ whenever $p,q\in G$ and $p\neq q$.
% \item[(2)] $\Delta(A_{pq})(1\otimes A_{q})=A_{p}\otimes A_{q}$ and $(A_{p}\otimes 1)\Delta(A_{pq})=A_{p}\otimes A_{q}$
% \end{enumerate}
% for all $p, q \in G$. Then $(A, \Delta)$ is called  a $G$-cograded multiplier Hopf algebra(see \cite{DVW05,YW11a}).
%
% We extend the Sweedler notation for a comultiplication in the following way:
% for any $p, q \in G$, $a\in A_{pq}$ and $a'\in A_{q}$, we write
% \begin{eqnarray*}
% \Delta_{p, q}(a)(1\otimes a') = a_{(1, p)} \otimes a_{(2, q)}a'.
% \end{eqnarray*}
%
% It seems that most of the theory of multiplier Hopf algebras can be generalized to the group-cograded case,
% such as integrals (see \cite{ADV07}), quasitriangular structure (see \cite{DVW05}).

\section{Drinfel'd twist}
\def\theequation{\thesection.\arabic{equation}}
\setcounter{equation}{0}

 In the paper \cite{AEG02}, the properties of Drinfel'd twist for finite-dimensional Hopf algebras were studied:
 given a Hopf algebra or quasitriangular Hopf algebra, one gets another such structure twisting it with a Drinfel'd twist,
 and the authors also determine how the integral of the dual to finite-dimensional unimodular Hopf algebra changes under a twist
 (see Theorem 3.4 in \cite{AEG02}).

 However, there is still a question here: for infinite-dimensional Hopf algebras, does Theorem 3.4 also holds?

 \subsection{Twist on multiplier Hopf algebras}

 Let $(A, \Delta)$ be a regular multiplier Hopf algebra, we first generalize the Drinfel'd twist to the multiplier Hopf algebra case,
 and then construct some new multiplier Hopf algebras by this twist.

 \textbf{Definition \thesection.1}
 A twist for regular multiplier Hopf algebra $A$ is an invertible element $J \in M(A \otimes A)$, which satisfies
 \begin{eqnarray}
 && (\Delta \otimes \iota)(J)(J\otimes 1) = (\iota \otimes \Delta)(J)(1 \otimes J). \label{0}
 \end{eqnarray}

 \textbf{Remark \thesection.2} (1) Apply $(\iota \otimes \varepsilon \otimes \iota)$ to the equation above, one sees that as in the Hopf case
 $c=(\varepsilon \otimes \iota)(J)=(\iota \otimes\varepsilon)(J)$ is a non-zero scalar for twist $J$.
 One can always replace $J$ by $c^{-1}J$ to normalize the twist in such way that
 \begin{eqnarray}
 && (\varepsilon \otimes \iota)(J)=(\iota \otimes\varepsilon)(J) = 1.
 \end{eqnarray}
 In the following, we will always assume that $J$ is normalized in this way.

 (2) Let $x \in M(A)$ be a invertible element such that $\varepsilon(x)=1$. If $J$ is a twist for $A$,
 then so is $J^{x}:= \Delta(x)J(x^{-1} \otimes x^{-1})$.
 Indeed, it is similar to the one in Hopf algebra case except that we should take the (unique) extension for the homomorphism
 $\iota \otimes \Delta$ and $\Delta \otimes \iota$ from $A \otimes A$ to $M(A \otimes A)$.
 The Twists $J$ and $J^{x}$ are said to be \emph{gauge equivalent}.
 \\

 \textbf{Example \thesection.3}
 Let $H$ be a Hopf algebra with a twist $J$ and $A$ be a multiplier Hopf algebra.
 Then $H\otimes A$ is a multiplier Hopf algebra with the product, coproduct, counit and antipode as follows.
 \begin{eqnarray*}
 && (h\otimes a)(h'\otimes a') = hh'\otimes aa', \quad \Delta(h\otimes a) = (\iota\otimes\tau\otimes\iota)(\Delta(h)\otimes \Delta(a)),\\
 && \varepsilon(h\otimes a)=\varepsilon(h)\varepsilon(a), \qquad S(h\otimes a) = S(h)\otimes S(a).
 \end{eqnarray*}
 In this case, there is a twist $J_{13}$ on $H\otimes A$,
 where $J_{13}=(\iota\otimes\tau\otimes\iota)(J\otimes 1_{M(A\otimes A)})$.
 \\

 Furthermore, suppose that $J(1 \otimes a),(1 \otimes a)J \in M(A)\otimes A$
 and $J(a \otimes 1), (a \otimes 1)J\in A\otimes M(A)$ for all $a\in A$,
 at this time, so is $J^{-1}$ and we call $J$ the general twisting for $A$.
 We denote $(a \otimes 1)J = a J^{(1)} \otimes J^{(2)}$ and $J(a \otimes 1) =J^{(1)}a \otimes J^{(2)}$ in $A \otimes M(A)$.

 In the above assumptions, we can define a multiplier $Q_{J}=S(J^{(1)})J^{(2)}\in M(A)$ by
 \begin{eqnarray*}
 && Q_{J}a=m(S \otimes \iota)(J(1\otimes a)) \in A, \\
 && aQ_{J}=m(S \otimes \iota)(J(S^{-1}(a) \otimes 1)) \in A
 \end{eqnarray*}
 for all $a\in A$. This multiplier $Q_{J}$ is invertible with the inverse $Q^{-1}_{J} =J^{-(1)} S(J^{-(2)})$. Indeed,
 for any $a, b \in A$
 \begin{eqnarray*}
  a Q^{-1}_{J}Q_{J} b
 &=& a J^{-(1)} S(J^{-(2)}) S(J^{(1)})J^{(2)} b \\
 &=& a J^{-(1)} S(J^{(1)}J^{-(2)})J^{(2)} b \\
 &=& a J'^{(1)} J^{-(1)} S(J^{(1)}J^{-(2)}) \varepsilon(J'^{(2)}) J^{(2)} b \\
 &=& a J'^{(1)} J^{-(1)} S(J'^{(2)}_{(1)} J^{(1)}J^{-(2)}) J'^{(2)}_{(2)} J^{(2)} b \\
 &=& a J^{(1)}_{(1)} J'^{(1)} J^{-(1)} S(J^{(1)}_{(2)} J'^{(2)} J^{-(2)}) J^{(2)} b \\
 &=& a J^{(1)}_{(1)} S(J^{(1)}_{(2)}) J^{(2)} b \\
 &=& ab
 \end{eqnarray*}
 and similarly $ a Q^{-1}_{J}Q_{J} b = ab$.

 The element $Q_{J}$ satisfies
 \begin{eqnarray}
 \Delta(Q_{J})=(S\otimes S)(J^{-1}_{21}) (Q_{J}\otimes Q_{J}) J^{-1}, \label{1}
 \end{eqnarray}
 where $J^{-1}_{21} = \tau J^{-1}$. These equations make sense because of the extension of (anti-) homomorphism
 (see Proposition A.5 in \cite{V94}).
 \\

 By the twist $J$, we can get a new multiplier Hopf algebra as follows, which generalizes the result in \cite{AEG02}.

 \textbf{Proposition \thesection.4}
 Let $(A, \Delta)$ be a regular multiplier Hopf algebra, then $(A^{J}, \Delta^{J})$ is also a regular multiplier Hopf algebra
 with the same algebra structure and counit as $(A, \Delta)$ and the comultiplication and antipode are given by
 \begin{eqnarray*}
 && \Delta^{J}(a) = J^{-1} \Delta(a) J, \\
 && S^{J}(a) = Q^{-1}_{J} S(a) Q_{J}
 \end{eqnarray*}
 for all $a\in A$.

 \emph{Proof} It is sufficient to check the equivalent definition, i.e., equations (\ref{1.1}) and (\ref{1.2}).
 Firstly, it is easy to check that $\Delta^{J}$ is a homomorphism and $S^{J}$ is a bijective anti-homomorphism.
 In the following, we only check (\ref{1.1}) and (\ref{1.2}) is similar.
 \begin{eqnarray*}
  (\varepsilon\otimes\iota)(\Delta^{J}(a)(1\otimes b))
 &=& (\varepsilon\otimes\iota) (J^{-1} \Delta(a) J (1\otimes b)) \\
 &=& (\varepsilon\otimes\iota) (J^{-1}) (\varepsilon\otimes\iota) (\Delta(a) (\varepsilon\otimes\iota) (J) b \\
 &=& ab,
 \end{eqnarray*}
 where the second equation holds because $\varepsilon$ ia a homomorphism.
 \begin{eqnarray*}
  m(S^{J}\otimes\iota)(\Delta^{J}(a)(1\otimes b))
 &=& m(S^{J}\otimes\iota)(J^{-1}\Delta(a)J(1\otimes b)) \\
 &=& m(S^{J}\otimes\iota)(J^{-(1)}a_{(1)}J'^{(1)} \otimes J^{-(2)}a_{(2)}J'^{(2)} b) \\
 &=& Q^{-1}_{J} S(J^{-(1)}a_{(1)}J'^{(1)}) Q_{J} J^{-(2)}a_{(2)}J'^{(2)} b \\
 &=& Q^{-1}_{J} S(J'^{(1)}) S(a_{(1)}) S(J^{-(1)}) Q_{J} J^{-(2)}a_{(2)}J'^{(2)} b \\
 &=& \varepsilon(a)b,
 \end{eqnarray*}
 where $J^{-1}\Delta(a)J(1\otimes b) = J^{-(1)}a_{(1)}J'^{(1)} \otimes J^{-(2)}a_{(2)}J'^{(2)} b$ makes sense,
 since $J(1\otimes b)\in A\otimes A$ and denote it as $J'^{(1)} \otimes J'^{(2)} b$,
 and then $\Delta(a)(J'^{(1)} \otimes J'^{(2)} b)\in A\otimes A$ by the "cover" technique shown in \cite{V08}.
 $\hfill \Box$
 \\

 In the subsequent subsection, we will use the following lemma, which is the direct deformation of equation (\ref{0}).

 \textbf{Lemma \thesection.5} Let $J$ be a twist for multiplier Hopf algebra $A$, and for any $a, b\in A$
 \begin{eqnarray*}
 && a S(J^{(1)})J^{(2)}_{(1)} \otimes b J^{(2)}_{(2)} = (a\otimes b)(Q_{J}\otimes 1)J^{-1}, \\
 && a J^{-(1)}S(J^{-(2)}_{(1)}) \otimes b S(J^{-(2)}_{(2)}) = (a Q^{-1}_{J} \otimes b) (S\otimes S)(J).
 \end{eqnarray*}

 \emph{Proof} We only check the first equation, and the second one is similar.
 \begin{eqnarray*}
  \big(a S(J^{(1)})J^{(2)}_{(1)} \otimes b J^{(2)}_{(2)} \big) J
 &=& a S(J^{(1)})J^{(2)}_{(1)} \bar{J}^{(1)} \otimes b J^{(2)}_{(2)} \bar{J}^{(2)} \\
 &\stackrel{(\ref{0})}{=}&  a S(J^{(1)}_{(1)} \bar{J}^{(1)} ) J^{(1)}_{(2)} \bar{J}^{(2)} \otimes b J^{(2)} \\
 &=& a S(\bar{J}^{(1)}) S(J^{(1)}_{(1)}) J^{(1)}_{(2)} \bar{J}^{(2)} \otimes b J^{(2)} \\
 &=& a S(\bar{J}^{(1)}) \bar{J}^{(2)} \otimes b \\
 &=& a Q_{J} \otimes b,
 \end{eqnarray*}
 so $\big(a S(J^{(1)})J^{(2)}_{(1)} \otimes b J^{(2)}_{(2)} \big) J = a Q_{J} \otimes b$,
 by multiplier $J^{-1}$ from the right side, we get the assertion.
 $\hfill \Box$

 \subsection{Quasitriangular structure under a twist}

 Recall from \cite{Z99}, a regular multiplier Hopf algebra $(A, \Delta)$ is called \emph{quasitriangular}
 if there exists an invertible multiplier $\mathcal{R}$ in $M(A\otimes A)$ which is subject to
 \begin{enumerate}
   \item[(1)] $(\Delta\otimes\iota)(\mathcal{R}) = \mathcal{R}^{13} \mathcal{R}^{23}$,
              $(\iota\otimes\Delta)(\mathcal{R}) = \mathcal{R}^{13} \mathcal{R}^{12}$,
   \item[(2)] $\mathcal{R}\Delta(a) = \Delta^{cop}(a)\mathcal{R}$ for all $a\in A$,
   \item[(3)] $(\iota\otimes\varepsilon)(\mathcal{R}) = 1 = (\varepsilon\otimes\iota)(\mathcal{R})$.
 \end{enumerate}

 Let $(A, \Delta)$ be a multiplier Hopf algebra, by Proposition \thesection.4, we get $(A^{J}, \Delta^{J})$ is also a multiplier Hopf algebra.
 When $(A, \Delta)$ is quasitriangular, how about $(A^{J}, \Delta^{J})$?
 \\

 \textbf{Theorem \thesection.6}
 Let $(A, \Delta)$ be a quasitriangular multiplier Hopf algebra with generalized $R$-matrix $\mathcal{R}$,
 then $(A^{J}, \Delta^{J})$ is also quasitriangular, and the quasitriangular structure given by
 \begin{eqnarray*}
 \mathcal{R}^{J} = J^{-1}_{21} \mathcal{R} J.
 \end{eqnarray*}

 \emph{Proof} It is sufficient to check the equation (1), (2) and (3) in the definition of quasitriangular.
 Firstly, we check $\mathcal{R}^{J}\Delta^{J}(a) = (\Delta^{J})^{cop}(a) \mathcal{R}^{J}$. In fact,
 \begin{eqnarray*}
  \mathcal{R}^{J}\Delta^{J}(a)
  &=& J^{-1}_{21} \mathcal{R} J J^{-1} \Delta(a) J = J^{-1}_{21} \mathcal{R} \Delta(a) J, \\
  (\Delta^{J})^{cop}(a) \mathcal{R}^{J}
  &=& \tau (J^{-1} \Delta(a) J) J^{-1}_{21} \mathcal{R} J
   = J^{-1}_{21} \Delta^{cop}(a) J_{21} J^{-1}_{21} \mathcal{R} J \\
  &=& J^{-1}_{21} \Delta^{cop}(a) \mathcal{R} J.
 \end{eqnarray*}
 Because $(A, \Delta)$ is quasitriangular, $\mathcal{R} \Delta(a) = \Delta^{cop}(a) \mathcal{R}$, then the equation holds.

 Secondly, by the extension of homomorphism $\varepsilon$, it is easy to check that
 $(\iota\otimes\varepsilon)(\mathcal{R}^{J}) = 1 = (\varepsilon\otimes\iota)(\mathcal{R}^{J})$.

 Finally, we need to check $(\Delta^{J}\otimes\iota)\mathcal{R}^{J} = {\mathcal{R}^{J}}^{13} {\mathcal{R}^{J}}^{23}$ and
 $(\iota\otimes\Delta^{J})\mathcal{R}^{J} = {\mathcal{R}^{J}}^{13} {\mathcal{R}^{J}}^{12}$.
 Here, we only show the proof of the first equation, and the second one is similar to verify.

 For any $x\otimes y \in A \otimes A$, $x\otimes y = \sum_{i}\Delta^{J}(a_{i})(1\otimes b_{i})$ because of Proposition \thesection.3.
 Thus
 \begin{eqnarray*}
 && (\Delta^{J}\otimes\iota)(\mathcal{R}^{J})(x\otimes y\otimes z) \\
 &=& (\Delta^{J}\otimes\iota)(\mathcal{R}^{J})(\sum_{i}\Delta^{J}(a_{i})(1\otimes b_{i}) \otimes z) \\
 &=& \sum_{i} (\Delta^{J}\otimes\iota)(\mathcal{R}^{J}(a_{i} \otimes z))(1\otimes b_{i} \otimes 1) \\
 &=& \sum_{i} (J^{-1}\otimes 1) (\Delta\otimes\iota)(J^{-1}_{21} \mathcal{R} J (a_{i} \otimes z)) (J\otimes 1) (1\otimes b_{i} \otimes 1) \\
 &=& \sum_{i} (J^{-1}\otimes 1) (\Delta\otimes\iota)(J^{-1}_{21}) (\Delta\otimes\iota) (\mathcal{R})
     (\Delta\otimes\iota)(J (a_{i} \otimes z)) (J\otimes 1) (1\otimes b_{i} \otimes 1) \\
 &=& \sum_{i} (J^{-1}\otimes 1) (\Delta\otimes\iota)(J^{-1}_{21}) \mathcal{R}^{13}\mathcal{R}^{23}
     (\Delta\otimes\iota)(J) (\Delta\otimes\iota) (a_{i} \otimes z) (J\otimes 1) (1\otimes b_{i} \otimes 1) \\
 &=& \sum_{i} (J^{-1}\otimes 1) (\Delta\otimes\iota)(J^{-1}_{21}) \mathcal{R}^{13}\mathcal{R}^{23} (\Delta\otimes\iota)(J)(J\otimes 1) \\
     &&  (J^{-1}\otimes 1) (\Delta\otimes\iota) (a_{i} \otimes z) (J\otimes 1) (1\otimes b_{i} \otimes 1) \\
 &=& \sum_{i} (J^{-1}\otimes 1) (\Delta\otimes\iota)(J^{-1}_{21}) \mathcal{R}^{13}\mathcal{R}^{23}
     (\iota\otimes\Delta)(J)(1 \otimes J)(J^{-1}\Delta(a_{i})J \otimes 1)(1\otimes b_{i} \otimes z) \\
 &=& (\dot{J}^{-(1)} J^{-(2)}_{(1)} \otimes \dot{J}^{-(2)} J^{-(2)}_{(2)} \otimes J^{-(1)}) \mathcal{R}^{13}\mathcal{R}^{23}
     (\iota\otimes\Delta)(J)(1 \otimes J)(x\otimes y \otimes z) \\
 &=& \underline{\dot{J}^{-(1)} J^{-(2)}_{(1)}} \mathcal{R}^{(1)} J'^{(1)} x
       \otimes \underline{\dot{J}^{-(2)} J^{-(2)}_{(2)}} \mathcal{R}'^{(1)} J'^{(2)}_{(1)} \dot{J}'^{(1)} y
       \otimes \underline{J^{-(1)}} \mathcal{R}^{(2)}\mathcal{R}'^{(2)} J'^{(2)}_{(2)} \dot{J}'^{(2)} z \\
 &=& J^{-(2)} \dot{J}^{-(1)}_{(2)} \mathcal{R}^{(1)} J'^{(1)} x
       \otimes \dot{J}^{-(2)}  \underline{\mathcal{R}'^{(1)} J'^{(2)}_{(1)}}  \dot{J}'^{(1)} y
       \otimes J^{-(1)} \dot{J}^{-(1)}_{(1)} \mathcal{R}^{(2)} \underline{\mathcal{R}'^{(2)} J'^{(2)}_{(2)}} \dot{J}'^{(2)} z \\
 &=& J^{-(2)} \underline{\dot{J}^{-(1)}_{(2)} \mathcal{R}^{(1)}}  J'^{(1)} x
       \otimes \dot{J}^{-(2)}  J'^{(2)}_{(2)}\mathcal{R}'^{(1)} \dot{J}'^{(1)} y
       \otimes J^{-(1)} \underline{\dot{J}^{-(1)}_{(1)} \mathcal{R}^{(2)}}  J'^{(2)}_{(1)}\mathcal{R}'^{(2)} \dot{J}'^{(2)} z \\
 &=& J^{-(2)} \mathcal{R}^{(1)} \underline{\dot{J}^{-(1)}_{(1)}  J'^{(1)}} x
       \otimes \underline{\dot{J}^{-(2)}  J'^{(2)}_{(2)}}  \mathcal{R}'^{(1)} \dot{J}'^{(1)} y
       \otimes J^{-(1)} \mathcal{R}^{(2)}  \underline{\dot{J}^{-(1)}_{(2)} J'^{(2)}_{(1)}}  \mathcal{R}'^{(2)} \dot{J}'^{(2)} z \\
 &=& J^{-(2)} \mathcal{R}^{(1)} J^{(1)} x
       \otimes \dot{J}^{-(2)} \mathcal{R}'^{(1)} \dot{J}'^{(1)} y
       \otimes J^{-(1)} \mathcal{R}^{(2)}  J^{(2)}  J'^{-(1)} \mathcal{R}'^{(2)} \dot{J}'^{(2)} z \\
 &=& {\mathcal{R}^{J}}^{13} {\mathcal{R}^{J}}^{23}(x\otimes y\otimes z),
 \end{eqnarray*}
 where the penultimate equation holds because of $(\Delta\otimes\iota)(J^{-1})(\iota\otimes\Delta)(J)=(J\otimes 1)(1\otimes J)$.
 $\hfill \Box$
 \\

 From Theorem \thesection.6 and Proposition 2.6 in \cite{DVW05}, we can easily get the following result.

 \textbf{Proposition \thesection.7}
 Let $(A, \Delta)$ be a quasitriangular multiplier Hopf algebra and $J$ a twist. Then for all $a\in A$
 \begin{eqnarray*}
 (S^{J})^{4}(a) = gag^{-1},
 \end{eqnarray*}
 where $g=\mu S(\mu)^{-1}$ and $\mu = S^{J}({\mathcal{R}^{J}}^{(2)}) {\mathcal{R}^{J}}^{(1)}$.

 \subsection{Integrals under a twist}

 Because the integral play a important role in the Pontryagin duality, there is naturally a question:
 for a regular multiplier Hopf algebra $(A, \Delta)$ with integrals, does $(A^{J}, \Delta^{J})$ also has a integral?
 \\

 Before answering this question, we first consider a multiplier $u_{J}=Q^{-1}_{J} S(Q_{J})$ in $M(A)$.
 By equation (\ref{1}), we have
 \begin{eqnarray*}
 \Delta(u_{J})
 &=& \Delta(Q^{-1}_{J} S(Q_{J})) = J (Q^{-1}_{J} S(Q_{J}) \otimes Q^{-1}_{J} S(Q_{J}))(S^{2}\otimes S^{2})(J^{-1}) \\
 &=& J (u_{J} \otimes u_{J})(S^{2}\otimes S^{2})(J^{-1}).
 \end{eqnarray*}

 \textbf{Theorem \thesection.8}
 Let $(A, \Delta)$ be a counimodular algebraic quantum group with a non-zero left (resp. right) integral $\varphi$ (resp. $\psi$)
 and $J$ be a twist. Then the elements $\varphi^{J} = u_{J} \rightharpoonup \varphi$ and $\psi^{J} = \psi \leftharpoonup u^{-1}_{J}$
 are non-zero left and right integrals on $(A^{J}, \Delta^{J})$ respectively.

 \emph{Proof}
 We need to check that $(\iota\otimes\varphi^{J})\Delta^{J}(a) = \varphi^{J}(a)1$, equivalently
 $S(\iota\otimes\varphi^{J})\Delta^{J}(a) = \varphi^{J}(a)1$. In deed, for any $x\in A$, there exist $b\in A$ such that $x=S(b)$, thus
 \begin{eqnarray*}
 && S(\iota\otimes\varphi^{J})(\Delta^{J}(a)) x \\
 &=& S(\iota\otimes u_{J} \rightharpoonup \varphi)( J^{-1} \Delta(a) J ) x \\
 &=& S(\iota\otimes \varphi)( J^{-1} \Delta(a) J (1\otimes u_{J}) ) S(b) \\
 &=& S(\iota\otimes \varphi)((b\otimes 1) J^{-1} \Delta(a) J (1\otimes u_{J}) ) \\
 &=& S(\iota\otimes \varphi)( b J^{-(1)} a_{(1)} \dot{J}^{(1)} \otimes J^{-(2)} a_{(2)} \dot{J}^{(2)} u_{J} ) \\
 &=& S(b J^{-(1)} a_{(1)} \dot{J}^{(1)}) \varphi(J^{-(2)} a_{(2)} \dot{J}^{(2)} u_{J}) \\
 &=& S(\dot{J}^{(1)}) \underline{S(a_{(1)})} S(J^{-(1)}) S(b) \underline{\varphi( a_{(2)} \dot{J}^{(2)} u_{J} S^{2}(J^{-(2)}) )} \\
 &\stackrel{(1)}{=}&
     S(\dot{J}^{(1)}) \dot{J}^{(2)}_{(1)} {u_{J}}_{(1)} S^{2}(J^{-(2)}_{(1)}) S(J^{-(1)}) S(b)
     \varphi( a \dot{J}^{(2)}_{(2)} {u_{J}}_{(2)} S^{2}(J^{-(2)}_{(2)}) ) \\
 &=& (\iota\otimes\varphi)\Big((1\otimes a) \big( \underline{S(\dot{J}^{(1)}) \dot{J}^{(2)}_{(1)} \otimes \dot{J}^{(2)}_{(2)}} \big)
     \Delta(u_{J}) \big( (S\otimes S) ( b \underline{\underline{J^{-(1)} S(J^{-(2)}_{(1)}) \otimes  S(J^{-(2)}_{(2)})}} ) \big)  \Big) \\
 &\stackrel{(2)}{=}&
     (\iota\otimes\varphi)\Big((1\otimes a) \big( (Q_{J} \otimes 1)J^{-1} \big)
     \Delta(u_{J}) \big( (S\otimes S) ( (b Q_{J^{-1}} \otimes 1) (S\otimes S)(J) ) \big)  \Big) \\
 &=& (\iota\otimes\varphi)\Big((1\otimes a) \big( (Q_{J} \otimes 1)J^{-1} \big)
     (J (u_{J} \otimes u_{J})(S^{2}\otimes S^{2})(J^{-1})) \\
    && \big((S^{2}\otimes S^{2})(J) (S(Q_{J^{-1}}) \otimes 1) (S(b)\otimes 1) \big)  \Big) \\
 &=& (\iota\otimes\varphi)\Big((1\otimes a) (Q_{J} \otimes 1)
     (u_{J} \otimes u_{J})(S(Q_{J^{-1}}) \otimes 1) (S(b)\otimes 1) \big)  \Big) \\
 &=& (\iota\otimes\varphi)\Big((1\otimes a) (1\otimes u_{J}) (S(b)\otimes 1) \big)  \Big) \\
 &=& \varphi^{J}(a)x,
 \end{eqnarray*}
 where equation (1) holds because of $S(\iota\otimes\varphi)(\Delta(a)(1\otimes b)) = (\iota\otimes\varphi)((1\otimes a)\Delta(b))$
 and (2) holds because of Lemma \thesection.5.
 $\hfill \Box$
 \\

\section*{Acknowledgements}

% The authors would like to thank the referee for his/her valuable comments.
 The work was partially supported by the National Natural Science Foundation of China (Grant No. 11226070),
 the Fundamental Research Fund for the Central Universities(Grant No.KYZ201322)
 and the NJAUF (No. LXY201201019, LXYQ201201103).

\end {document}